\documentclass[12pt]{amsart}

\usepackage{amsmath,amsfonts,stackrel,amsthm,enumerate,subcaption}
\usepackage{graphicx,setspace,mathrsfs,placeins,mathtools,amssymb}
\usepackage[margin=1in]{geometry}
\usepackage[mathscr]{eucal}
\theoremstyle{definition}
\usepackage[american]{babel}
\usepackage{tikz}
\usetikzlibrary{calc,decorations.markings}
\usepackage{csquotes}
\allowdisplaybreaks

\numberwithin{equation}{section}
\newcommand{\ncom}{\newcommand}
\ncom{\nno}{\nonumber}
\ncom{\vone}{\vskip 2ex}
\ncom{\norm}{\|\;\;\|}
\ncom{\vspan}[1]{{{\rm\,span}\{ #1 \}}}
\ncom{\dm}[1]{ {\displaystyle{#1} } }
\ncom{\ri}[1]{{#1} \index{#1}}

\newtheorem{theorem}{\bf Theorem}[section]
\newtheorem{remark}{\bf Remark}[section]

\newtheorem{definition}{Definition}[section]
\newtheoremstyle
{remarkstyle}
{}
{11pt}
{}
{}
{\bfseries}
{:}
{     }
{\thmname{#1} \thmnumber{#2} }

\theoremstyle{remarkstyle}

\makeatletter
\newcommand\fs@norules{\def\@fs@cfont{\bfseries}\let\@fs@capt\floatc@ruled
	\def\@fs@pre{}%
	\def\@fs@post{}%
	\def\@fs@mid{\kern3pt}%
	\let\@fs@iftopcapt\iftrue}
\makeatother
\usepackage[ruled, lined, linesnumbered, commentsnumbered, resetcount, longend]{algorithm2e}
\usepackage{caption}
\usepackage{subcaption}

\begin{document}
\title[]{Parameter estimation of incomplete gamma subordinators   }
\author[]{Meena Sanjay Babulal$^*$$^a$, Sunil Kumar Gauttam$^*$$^b$ and Aditya Maheshwari$^\#$$^c$}
\email{$^a$19pmt002@lnmiit.ac.in, $^b$sgauttam@lnmiit.ac.in, $^c$adityam@iimidr.ac.in}

\address[]{
	$^*$Department of Mathematics, The LNM Institute of Information Technology, Rupa ki Nangal, Post-Sumel, Via-Jamdoli
	Jaipur 302031,
	Rajasthan, India. }
\address[]{$^\#$Operations Management and Quantitative Techniques Area,
	Indian Institute of Management Indore, Indore 453556, Madhya Pradesh, India.}
\begin{abstract}
In this paper, we estimate the parameters of InG, InG-$\epsilon$ and TInG subordinators which have been studied by Babulal \textit{et al} (see \cite{babulal}).  We have modified the  method of moments technique to use fractional moments  of 
the InG and InG-$\epsilon$ subordinator due to their infinite moments. 
For the TInG subordinator's parameter estimation, we have used the method of moments. We also compute the maximum likelihood estimator(MLE) for the parameter $\alpha$ of the InG and InG-$\epsilon$ subordinators using jump distribution of the process. We also discussed the asymptotic normality of MLE. 
\end{abstract}
  \maketitle

	\section{Introduction}

Recently, Beghin and Ricciuti (see \cite{ricciuti}) introduced incomplete gamma subordinators using  the lower-incomplete gamma function and represented it  as compound Poisson process with positive jumps. 
A tempered version TInG subordinator is also considered in order to overcome the drawback of infinite moments of InG and InG-$\epsilon$ subordinators. We (see \cite{babulal})  have subordinated general L\' evy process by these subordinators and studied their's mean, variance, correlation and other distributional properties along with long range dependence (LRD) property. \\

 In statistical analysis, parameter estimation plays a key role as it deals with calculating the values of parameters inside a statistical model using observable data.  
There are models which have real life application but they  depends on parameter which can not be measured directly. In such cases parameter estimation play vital role in making model applicable. \\

There are several other methods for parameter estimation like  Bayesian estimation (see \cite{ALJEDDANI2023221,Murphy2015BayesianPE}), 
Machine learning (ML)-based parameter estimation (see \cite{Bhuvaneswari2019MachineLP}), Markov chain Monte Carlo method (see \cite{Chalana1997ParameterEI}) and  Cramer-von–
Mises (CVM) estimation (see \cite{10.1111/j.2517-6161.1971.tb00884.x}) etc. Each method has its own strengths and limitations. Cahoy \textit{et
	al.} (see \cite{Cahoy2010}) have used method of moment technique to estimate parameters of the fractional Poisson process using first two moments of its interarrival times. The parameter estimation of fractional Yule process have been studied by Cahoy and Polito (see \cite{Cahoy2012SimulationAE})  using  method of moments technique on  first two moments of their sojourn time. Cahoy and Polito (see \cite{articlefbpdpPE}) used sojourn time of the fractional birth and the fractional death process to estimate theirs parameters using linear regression method. The method of moments technique is used by Stephan Haug \textit{et
	al.} (see \cite{garch})  on equidistance observation of the continuous time GARCH(1,1) process  to estimate its parameters using  first two moments and autocorrelation function of equidistance observation. 
    Andrew W. (see \cite{MLEitoP}) have used MLE to estimate parameters of the generalized Itô processes with discretely sampled data.\\
    
In this paper, we estimate the parameters of InG, InG-$\epsilon$ and TInG subordinator. The mean and variance of the InG and InG-$\epsilon$ subordinator doesn't exist. Therefore we have  modified the method of moments technique to use asymptotic behaviour of the fractional moments of InG and InG-$\epsilon$ subordinator to estimate their  parameter. We also compute maximum likelihood estimator (MLE) for parameter of InG and InG-$\epsilon$ subordinator. The asymptotic normality of the MLE is also discussed. We have used method of moment to estimate the parameter of TInG subordinator using it's first two moments. \\

Next, we study asymptotic normality 
 for the MLE of InG and InG-$\epsilon$ subordinator's parameter. 
The idea of asymptotic normality is important  because it demonstrates that, under certain conditions, estimator's distributions will converge to a normal distribution as the sample size approaches infinite. This property facilitates the development of confidence intervals and hypothesis tests, allowing statisticians to make robust inferences about population parameters. Understanding the theoretical underpinnings of parameter estimation and its asymptotic properties is essential for effective statistical practices.\\

The paper is organized as follows. In Section $2$, we have given some notation, definitions, and results which will be used later in the paper. In Section $3$, we compute the maximum likelihood estimator(MLE) for the parameter $\alpha$ of
the InG and InG-$\epsilon$ subordinators. In Section $4$, we studied the asymptotic normality for the estimator of the InG and InG-$\epsilon$ subordinator's parameter. The Section $5$ deals with parameter estimation of InG, InG-$\epsilon$, and TInG subordinators using the method of moments. 
\section{preliminaries}
	\noindent 
\setcounter{equation}{0}
	In this section, we introduce the notation, definition, and results that will be used later.
	Let $\mathbb{Z}^{+} :=\left\{0,1, 2,\cdots\right\}$ be the set of non-negative integers. Let $\left\{N (t, \lambda)\right\}_{t\geq0}$ be a Poisson process with rate $\lambda > 0$, so that
	\begin{align}
		p(n|t, \lambda):=\mathbb{P}[N (t, \lambda)=n]=\frac{(\lambda t)^{n}e^{-\lambda t}}{n!}, \qquad n \in \mathbb{Z}^{+}.
	\end{align}
	For simplicity of notation we write $N( t,\lambda)$ as $N(t)$, when no confusion arises.
	For $\alpha \in (0, 1]$, the  InG subordinator $\lbrace S_{\alpha}(t)\rbrace_{t \geq 0} $   (see \cite{ricciuti}) can be represented as a compound Poisson process
	\begin{align}\label{subordinator InG}
		S_{\alpha}(t) =
		\sum_{j=1}^{N_{\alpha}(t)}	Z_{j}^{\alpha},
	\end{align}
	where $\left\lbrace N_{\alpha}(t)\right\rbrace_{ t \geq 0  }$
	is a homogeneous Poisson process with the rate $\lambda := \alpha
	\Gamma (\alpha) $ and the jumps $Z_{j}^{\alpha}$ are i.i.d. random variables, taking values in $[1, +\infty)$, with probability
	density function
	\begin{align}\label{Zpdf}
		f_{Z^{\alpha}}(z) = \frac{(z-1)^{-\alpha}z^{-1}1_{z\geq1}}{\Gamma (1-\alpha)\Gamma(\alpha)} = \frac{\sin(\pi \alpha)1_{z\geq1} }{\pi (z-1)^{\alpha}z },   \quad   \alpha \in (0,1).
	\end{align}

	When $\alpha = 1$, the jumps are unitary, and the process coincides with the Poisson process (see \cite{ricciuti}).	
	Note that the subordinator $\lbrace S_{\alpha}(t)\rbrace_{t \geq 0} $ 
	have jumps of size greater than or equal to $1$.\\
	Similarly, the InG-$\epsilon$ subordinator $\lbrace S_{\alpha}^{(\epsilon)}(t)\rbrace_{t \geq 0}$ can be represented as a compound
	Poisson process
	\begin{equation}\label{subordinator InGE}
		S_{\alpha}^{(\epsilon)}(t) =\sum_{j=1}^{N^{\epsilon}(t)}	Z_{j}^{(\alpha,\epsilon)},
	\end{equation}
	where $N^{\epsilon} := \lbrace N^{\epsilon}(t)\rbrace_{ t \geq 0 }$ 
	is a homogeneous Poisson process with the rate $\lambda := \alpha
	\Gamma (\alpha) \epsilon^{-\alpha} $ and the jumps $Z_{j}^{(\alpha,\epsilon)}$ are i.i.d. random variables, taking values in $[\epsilon, +\infty)$, with probability density function
	\begin{align}\label{ZE pdf}
		f_{Z_{j}^{(\alpha,\epsilon)}}(z) = \frac{\epsilon^{\alpha}(z-\epsilon)^{-\alpha}z^{-1}1_{z\geq\epsilon}}{\Gamma (1-\alpha)\Gamma(\alpha)} 
		,   \quad   \alpha \in (0,1).
	\end{align}
	In contrast to the InG subordinator, the InG-$\epsilon$ subordinator $\lbrace S_{\alpha}^{(\epsilon)}(t)\rbrace_{t \geq 0}$ have jumps of size greater than or equal to $\epsilon$ (see \cite{ricciuti}). \\
	The TInG subordinator $\lbrace S_{\alpha, \theta}(t)\rbrace_{t \geq 0}$ can be represented as a compound
	Poisson process

	\begin{equation}\label{Tempered subordinator}
		S_{\alpha, \theta}(t)=	\sum_{j=1}^{N_{ \alpha, \theta}(t)}	Z_{j}^{ \alpha, \theta}, 
	\end{equation}
	where $ N_{\alpha, \theta}$ := $\left\lbrace N_{\alpha, \theta}(t)\right\rbrace_{t\geq 0} $ is a homogeneous Poisson process with rate $\lambda := \alpha \Gamma(\alpha;\theta)$, where $\Gamma(\alpha;\theta)$ is the upper incomplete gamma function defined as
	\begin{equation}\label{InGfU}
		\Gamma(\alpha;\eta)=\int_\eta^\infty e^{-y}y^{\alpha-1}dy,~~ \eta>0,0<\alpha\leq 1.\nonumber
	\end{equation}
	
    The jumps $Z_{j}^{ \alpha, \theta}$
	are i.i.d. random variables, taking values in $[1, +\infty)$ and with the probability density function
	\begin{equation}\label{pdf of ztheta}
		f_{Z_{j}^{ \alpha, \theta}}=\frac{e^{-\theta z}(z-1)^{-\alpha}z^{-1}1_{z\geq1}}{\Gamma (1-\alpha)\Gamma(\alpha;\theta)}, \quad \alpha \in (0,1).
	\end{equation} 
	Observe that
	the mean for subordinators the InG and InG-$\epsilon$ does not exist, but  mean and variance of the TInG  subordinator $S_{\alpha, \theta}(t)$ are given by (see \cite{ricciuti})
	\begin{equation}\label{mean tempered subordinator}
		\mathbb{E} S_{\alpha, \theta}(t) = t\alpha\theta^{\alpha-1}e^{-\theta}, \end{equation}
	\begin{equation}\label{var tempered subordinator }
		\mbox{Var}  S_{\alpha, \theta}(t)= t\alpha\theta^{{\alpha-1}}e^{-\theta}+t({\alpha-1})\alpha\theta^{\alpha-2}e^{-\theta}.
	\end{equation}

\section{ Maximum likelihood estimation }
	In this section, we compute the maximum likelihood estimator(MLE) for the parameter $\alpha$ of the InG and InG-$\epsilon$ subordinators. 

    Our likelihood functions are based on the random samples from pdfs \eqref{Zpdf} and \eqref{ZE pdf}. First, we illustrate the  details for computation of MLE for the random sample from the density \eqref{Zpdf}.

	Suppose $Z_1, Z_2, ..., Z_n$ be iid observation on density (\ref{Zpdf}) and $z_1, z_2, ..., z_n$ are realised value of random sample $Z_1, Z_2, ..., Z_n$. 
  
    For computational ease, maximization is performed by optimizing the log-likelihood function  $\ln L(\alpha)$, which transforms the product of probabilities into a sum and facilitates computation while preserving the location of the optimum. 
	
	\begin{align}
		\ln L(\alpha) = \sum_{i=1}^{n} \ln \frac{\sin(\pi \alpha)1_{z_i\geq1} }{\pi (z_i-1)^{\alpha}z_i }
	\end{align}
	
	To find possible candidates for MLE, we take derivative of $\ln L(\alpha)$ with respect to $\alpha$ and equate it equal to $0.$
	\begin{align*}
		\frac{d}{d \alpha} \ln L(\alpha) &= \frac{d}{d \alpha} \left[\sum_{i=1}^{n} \ln \frac{\sin(\pi \alpha)1_{z_i\geq1} }{\pi (z_i-1)^{\alpha}z_i }\right] \\
		&= \sum_{i=1}^{n} \frac{d}{d \alpha} \ln \left[\frac{\sin(\pi \alpha)1_{z_i\geq1}}{\pi (z_i-1)^{\alpha}z_i}\right] \\
		&= \sum_{i=1}^{n} \left[\frac{\pi \cos(\pi \alpha)}{\sin(\pi \alpha)} - \ln(z_i-1) 
		\right] \\
		&= \sum_{i=1}^{n} \left[\frac{\pi}{\tan(\pi \alpha)} - \ln(z_i-1) 
		\right]
	\end{align*}
	

\begin{align}\label{InG equation of critical point}
\frac{d}{d \alpha} \ln L(\alpha) &=0 \nonumber\\
 \frac{n\pi}{\ln(\prod_{i=1}^{n} (z_i-1))}& = {\tan(\pi \alpha)}
\end{align}




We want to find the values of $\alpha$, which satisfies equation \eqref{InG equation of critical point}. 
Whenever arctan takes value is betweeen $(0,\pi/2]$, we get $\alpha$ value between $(0,1/2)$ and for value between $[-\pi/2,0)$ we get $\alpha$ value between $(1/2,1)$. Therefore, the unique solution $\alpha_0$ of the likelihood equation \eqref{InG equation of critical point}  is given by
\begin{equation}
	\alpha_0  = \left\{ \begin{array}{ll}
		
		\frac{1}{\pi} \arctan \left( \frac{n\pi}{\ln(\prod_{i=1}^{n} (z_i-1))} \right) & \mbox{ if $0< \arctan\left( \frac{n\pi}{\ln(\prod_{i=1}^{n} (z_i-1))} \right)\leq \dfrac{\pi}{2}$}\\
		\frac{1}{\pi} \left[\pi+\arctan \left( \frac{n\pi}{\ln(\prod_{i=1}^{n} (z_i-1))} \right)\right] & \mbox{if $-\dfrac{\pi}{2}\leq \arctan\left( \frac{n\pi}{\ln(\prod_{i=1}^{n} (z_i-1))} \right)\leq 0.$}
	\end{array}\right.\\
\end{equation}

The second derivative of $\ln L(\alpha)$ is $-n\pi^{2}\csc^2(\pi \alpha)$, which is negative for each $\alpha$ in range $(0,1).$
 Hence ${\alpha}_{0}$ is a local maxima.\\
 
 Note that $\alpha_0$ is the only critical point in the interior of the parameter space $(0,1)$ and a local maxima occur on it.
 For boundary points  $\alpha=0$ and $\alpha=1$ the function $L(\alpha)$ is $0$. Therefore it is the point of global maximum of $L(\alpha)$ over $(0,1)$. So we have the following theorem.

\begin{theorem}
    The maximum likelihood estimator of parameter $\alpha$ for InG subordinator is given by 
    \begin{equation}\label{InG mle}
	\hat{\alpha}_{mle}  = \left\{ \begin{array}{ll}
		
		\frac{1}{\pi} \arctan \left( \frac{n\pi}{\ln(\prod_{i=1}^{n} (Z_i-1))} \right) & \mbox{$0< \arctan\left( \frac{n\pi}{\ln(\prod_{i=1}^{n} (Z_i-1))} \right)\leq \dfrac{\pi}{2}$}\\
		\frac{1}{\pi} \left[\pi+\arctan \left( \frac{n\pi}{\ln(\prod_{i=1}^{n} (Z_i-1))} \right)\right] & \mbox{$-\dfrac{\pi}{2}\leq \arctan\left( \frac{n\pi}{\ln(\prod_{i=1}^{n} (Z_i-1))} \right)\leq 0$}
	\end{array}\right.
\end{equation}
where $Z_1,\cdots, Z_n$ is a random sample from the pdf \eqref{Zpdf}.
\end{theorem}

On similar line of argument, we have the following result regarding the MLE of parameter $\alpha$ of InG-$\epsilon$ subordinator. 
 
\begin{theorem}
    The MLE of the paramtere $\alpha$  for InG-$\epsilon$ subordinator is given by 

\begin{equation}\label{InG-e mle}
\hat{\alpha}^{(\epsilon)}_{mle}  = \left\{ \begin{array}{ll}
	
	\frac{1}{\pi} \arctan \left( \frac{n\pi}{\ln(\prod_{i=1}^{n} (Z_i-\epsilon))} \right) & \mbox{if $0< \arctan \left( \frac{n\pi}{\ln(\prod_{i=1}^{n} (Z_i-\epsilon))} \right)\leq \dfrac{\pi}{2}$}\\
		\frac{1}{\pi} \left[\pi+\arctan \left( \frac{n\pi}{\ln(\prod_{i=1}^{n} (Z_i-\epsilon))} \right)\right] & \mbox{if $\dfrac{-\pi}{2}\leq \arctan \left( \frac{n\pi}{\ln(\prod_{i=1}^{n} (Z_i-\epsilon))} \right)\leq 0$}
\end{array}\right.
\end{equation}
	where $Z_1,\cdots, Z_n$ is a random sampple from the pdf \eqref{ZE pdf}.\\ 
\end{theorem}

\section{Asymptotic normality of MLE}

In this section, we will discuss the asymptotic normality of the estimator \eqref{InG mle} and \eqref{InG-e mle} 
using the Cramer theorem (see \cite{inbook}). First, we state  definition of the asymptotic normality (see \cite{inbook}).\\
\begin{definition}
    
A sequence of estimators  ($T_n$) is said to be asymptotically normal with mean $\tau (\theta)$  and variance $\frac{\nu(\theta)}{n}$, and write $T_n$ is  AN$(\tau (\theta),\frac{\nu(\theta)}{n})$ for $\theta \in \mathcal{O}$, where $\mathcal{O}$ is any open interval of real number $\mathbb{R}$. \textit{i.e.} \
    $$\dfrac{T_n-\tau (\theta)}{\sqrt{(\nu(\theta)/n)}}\overset{d}{\rightarrow} \mathcal{N}(0,1), \qquad \mbox{as } n\rightarrow \infty.$$
\end{definition}

We use following Cramer's theorem to verify asymptotic normality of MLE 
(see \cite{inbook}). 

 \begin{theorem}     
     Let   
    $ \left\{     f_\theta | \theta \in \mathcal{O}\right\}$ is a family of pdfs, where $\mathcal{O} $   is any open interval of real number $\mathbb{R}$. Then \\

\begin{enumerate}

 \item[(i)] \( \dfrac{\partial \log f_\theta}{\partial \theta}, \, \dfrac{\partial^2 \log f_\theta}{\partial \theta^2}, \, \dfrac{\partial^3 \log f_\theta}{\partial \theta^3} \) exist for all \( \theta \in \Theta \). Also,\\
    \[
     \int_{-\infty}^\infty \frac{\partial \log f_\theta}{\partial \theta} \,  =   \mathbb{E}_\theta \left[ \frac{\partial \log f_\theta(X)}{\partial \theta} \right] = 0 \quad \text{for all } \theta \in \Theta.\\
    \]

    \item[(ii)] \(\int_{-\infty}^\infty \frac{\partial^2 \log f_\theta}{\partial \theta^2}  \, dx = 0 \) for all \( \theta \in \Theta. \)\\

    \item[(iii)] \(  -\infty<\int_0^\infty \frac{\partial^2 \log f_\theta}{\partial \theta^2} f_\theta(x)\, dx < 0 \) for all \( \theta \in \Theta. \)\\

    \item[(iv)] There exists a function \( H(x) \) such that
    \[
    \left| \frac{\partial^3 \log f_\theta(x)}{\partial \theta^3} \right| < H(x) \quad \text{and} \quad \int_{-\infty}^\infty H(x) f_\theta(x) dx < \infty.
    \]

\end{enumerate}

Conditions (i) - (iv) imply that a consistent solution \( \hat{\theta}_n \) of the likelihood equation is asymptotically normal, that is,
\[
\sigma^{-1}\sqrt{n}(\hat{\theta}_n - \theta_0) \xrightarrow{d} Z,
\]
where \( Z \sim N(0, 1) \), and
\[
\sigma^2 = \mathbb{E}_\theta \left[ \left( \frac{\partial \log f(x)}{\partial \theta} \right)^{2} \right]^{-1}.\\
\]
 \end{theorem}

We verify the conditions $(i)-(iv)$ of theorem for the MLE \eqref{InG mle}.

\hspace{1cm}(i)
\begin{align}\label{cond 1}
\dfrac{\partial \log f_{Z_\alpha}(z)}{\partial \alpha}&=\pi \cot(\alpha \pi) - \log(-1 + z),\\
 \dfrac{\partial^{2} (\log f_{Z_\alpha}(z))}{\partial \alpha^{2}}&=-\pi^2 \csc^2(\alpha \pi),\nonumber\\
   \frac{\partial^{3} (\log f_{Z_\alpha}(z))}{\partial \alpha^{3}}&=2 \pi^3 \cot(\alpha \pi) \csc^2(\alpha \pi).\nonumber
\end{align}

To show $\mathbb{E}\left[ \dfrac{\partial}{\partial \alpha} \log f_{Z_\alpha}(z) \right]=0$, from equation \eqref{cond 1} it is enough to prove that mean of \(\log(Z - 1)\) is equal to $\pi \cot(\pi \alpha)$. 
The mean of \(\log(Z - 1)\) is given by
\[
E[\log(Z - 1)] = \int_{1}^{\infty} \log(z - 1) f_{Z_\alpha}(z) \, dz.
\]
Using the expression \eqref{Zpdf} for the density $f_{Z_\alpha}(z)$, we get 
\[
\mathbb{E}[\log(Z - 1)] = \int_{1}^{\infty} \log(z - 1) \left( \frac{\sin(\pi \alpha)}{\pi} \cdot \frac{(z - 1)^{-\alpha}}{z} \right) dz.
\]
\if
Simplifying the integral, this expression simplifies to
\[
\mathbb{E}[\log(Z - 1)] = \frac{\sin(\pi \alpha)}{\pi} \int_{1}^{\infty} \log(z - 1) \frac{(z - 1)^{-\alpha}}{z} \, dz.\]
\fi
By change of variable, we get
\[
\mathbb{E}[\log(Z - 1)] = \frac{\sin(\pi \alpha)}{\pi} \int_{0}^{\infty} \frac{\log(u) 
 u^{-\alpha}}{u + 1} \, du
=\pi cot[\alpha \pi]. \]\\

This completes the verification of condition $(i)$.\\
\begin{align*}
    &(ii) \qquad \frac{\partial^{2} f_{Z_\alpha}(z)}{\partial \alpha^{2}}=\frac{-(2 (-1 + z)^{-\alpha} \cos[\alpha\pi] \log[-1 + z])}{z}
     -\frac{(\pi (-1 + z)^{-a} \sin[
   \alpha \pi])}{z}\\
   &\hspace{4cm}+ \frac{((-1 + z)^{-\alpha} \log[-1 + z]^2 \sin[\alpha \pi])}{\pi z}.\\\\
   &\mbox{ To verify condition $(ii)$, we integrate second derivative of density function and get }\\\\
   &\qquad \qquad  \int_1^\infty  \frac{\partial^{2} f_{Z_\alpha}(z)}{\partial \alpha^{2}}dz=0.\\\\
   &(iii)  \qquad\int_1^\infty \frac{\partial^2 \log f_{Z_\alpha}(z)}{\partial \alpha^2} f_{Z_\alpha}(z)\, dz=-\pi^2 csc[\alpha \pi]^2<0.\\\\
   &(iv) \qquad \frac{\partial^{3}(\log f_{Z_\alpha}(z))}{\partial \alpha^{3}} =2 \pi^3 \cot(\alpha \pi) \csc^2(\alpha \pi).\\
\end{align*}
We choose function 
      $ H(z)=2 \pi^3 \cot(\alpha \pi) \csc^2(\alpha \pi)+1 $
      such that
   $ 
    \left| \dfrac{\partial^3 \log f_{Z_\alpha}(z)}{\partial \alpha^3} \right| < H(z) $
   and 
    $$\quad \int_{-\infty}^\infty H(z) f_{Z_\alpha}(z) dz =2 \pi^3 \cot(\alpha \pi) \csc^2(\alpha \pi)+1< \infty.$$
    Hence by Cramer's theorem, we have
   $$ 
\sigma^{-1}\sqrt{n}(\hat{\alpha}_{mle,n} - \alpha) \xrightarrow{d} Z $$
where $ Z \sim N(0, 1) $, and
$$
\sigma^2 = \mathbb{E}_\alpha \left[ \left( \frac{\partial \log f_{Z_\alpha}(z)}{\partial \alpha} \right)^{2} \right]^{-1}= \frac{2 \sin({\alpha \pi})}{\pi^{2} (3+\cos{(2\alpha \pi)}) }.$$
\ifx
\begin{remark}

    On similar line of arguments, global maximum of likelihood function $$L(\alpha) = \prod_{i=1}^{n} f( z_{i}^{(\alpha,\epsilon)}|\alpha)$$ occurs at a unique point $\alpha_0$, where 

\begin{equation*}
\alpha_0  = \left\{ \begin{array}{ll}
	
	\frac{1}{\pi} \arctan \left( \frac{n\pi}{\ln(\prod_{i=1}^{n} (z_i-\epsilon))} \right) & \mbox{$0\leq \dfrac{n\pi}{\ln(\prod_{i=1}^{n} (z_i-\epsilon))}\leq \dfrac{\pi}{2}$}\\
		\frac{1}{\pi} \left[\pi+\arctan \left( \frac{n\pi}{\ln(\prod_{i=1}^{n} (z_i-\epsilon))} \right)\right] & \mbox{$\dfrac{-\pi}{2}\leq \dfrac{n\pi}{\ln(\prod_{i=1}^{n} (z_i-\epsilon))}\leq 0.$}
\end{array}\right.
\end{equation*}
\end{remark}
\fi
\begin{remark}
     On similar line  using Cramer's theorem, asymptotic normality of estimator $\hat{\alpha}^{(\epsilon)}_{mle}$  is
   $$ 
\sigma^{-1}\sqrt{n}(\hat{\alpha}^{(\epsilon)}_{mle,n} - \alpha) \xrightarrow{d} Z $$

where $ Z \sim N(0, 1) $, and
$$\sigma^2 = \frac{1}{\left( \pi^2 + 2 \pi^2 \cot(\alpha \pi)^2 + 2 \pi \cot(\alpha \pi) \log(\epsilon) + \log(\epsilon)^2 \right) }.\\\\$$
\end{remark}

\section{Method of moments Estimation}

 In this section, we compute the MoM estimates for parameters of InG, InG-$\epsilon$ and TInG subordinators. The MoM is a straightforward approach used in parameter estimation, where we equate theoretical moments to their sample moments to estimate parameters. \\

 First, we compute MoM estimates for the parameters 
 of TInG subordinator using its mean and variance. The mean and variance of the TInG  subordinator $S_{\alpha, \theta}(t)$ are given by equations \eqref{mean tempered subordinator} and \eqref{var tempered subordinator }, respectively.


Let $X_{1}, X_{2},\cdots, X_{n}$ be sample from TInG subordinator  and we equate first and second sample moments with the subordinator's moments.

\begin{equation}\label{first moment tempered subordinator }
   \overline{X}=\frac{1}{n}\sum_{i=1}^{n} X_{i}=t\alpha\theta^{\alpha-1}e^{-\theta}\end{equation}
   \begin{equation}\label{second moment tempered subordinator }
    \frac{1}{n}\sum_{i=1}^{n} X_{i}^{2}=t\alpha\theta^{{\alpha-1}}e^{-\theta}+t({\alpha-1})\alpha\theta^{\alpha-2}e^{-\theta}+(t\alpha\theta^{\alpha-1}e^{-\theta})^2.\end{equation}
From equations \eqref{first moment tempered subordinator } and \eqref{second moment tempered subordinator }, we can easily see that in both equations $\alpha$ and $\theta$  are inseparable and the equations have a very complex structure which makes them difficult to solve analytically. 
  We solve the above equation numerically and following are the details of our experiment.\\\\
To find estimator $\hat{\alpha}$ and $\hat{\theta}$ for the parameters of TInG subordinator, we generate $n=100$ samples of TInG subordinator using  its algorithm (see \cite{babulal}).  Then we calculate $\hat{\alpha}$ and $\hat{\theta}$ estimates for the  sample size $N=100$, $500$ and $1000$ and we take  average of their estimates to compute  
$\hat{\alpha}_{avg}$ and $\hat{\theta}_{avg}$.  Here we have used matlab inbuilt function to estimate parameters. These estimated values are shown along with their mean absolute deviation(MAD) and mean squared error(MSE) in Table  $1$-$4$.\\\\

It is easy to observe from Table $1$-$4$ that as sample size increases parameter estimator keep approaching true value. Here, we must note that sample size taken is small as compare to  real-world applications, including network traffic data where sizes ranges in millions. \\\\
Table $1$
Mean estimates of and dispersions from the true parameter for a simulated TInG subordinator data with ($\alpha$, $\theta$)=($0.1$, $0.2$).\\

\begin{tabular}{c c c c}
  &$\dfrac{N=100}{Mean \quad MAD  \quad {MSE}}$ &  $\dfrac{N=500}{Mean \quad MAD \quad {MSE}}$ &  $\dfrac{N=1000}{Mean \quad MAD \quad {MSE}}$
   \vspace{12pt}\\
  
  \hspace{5pt}$\hat{\alpha}$ & 
  0.1325\quad0.0109\quad0.000177 &
  0.1312\quad0.0093\quad0.0001321
  & 0.1313\quad0.0096\quad0.0001436 
 \vspace{12pt}\\

 \vspace{5pt} $\hat{\theta}$ & 0.2026\quad0.0120\quad0.000214 & 0.2018\quad0.0104\quad0.000166 & 0.2020\quad0.0106\quad0.00017\\
\end{tabular}
\\\\

Table $2$
Mean estimates of and dispersions from the true parameter for a simulated TInG subordinator data with ($\alpha$, $\theta$)=($0.5$, $0.4$) .\\\\

\begin{tabular}{c c c c}
  &$\dfrac{N=100}{Mean \quad MAD \quad {MSE}}$ &  $\dfrac{N=500}{Mean \quad MAD \quad {MSE}}$ &  $\dfrac{N=1000}{Mean \quad MAD \quad {MSE}}$\vspace{12pt}\\
  $\hat{\alpha}$ &  0.5730\quad0.0298\quad0.0013 & 0.5710 \quad 0.0283 \quad 0.0012 &

  0.5691\quad0.0291\quad0.0013
 \vspace{12pt}\\
  $\hat{\theta}$ & 0.3846\quad 0.0114\quad 0.000198  & 0.3843 \quad 0.0108 \quad 0.00018 &  0.3846\quad 0.0113\quad 0.00029 
  \vspace{12pt}\\
\end{tabular}\\\\

Table $3$
Mean estimates of and dispersions from the true parameter for a simulated TInG subordinator data with ($\alpha$, $\theta$)=($0.7$, $0.5$).\\\\

\begin{tabular}{c c c c}
  &$\dfrac{N=100}{Mean \quad MAD  \quad {MSE}}$ &  $\dfrac{N=500}{Mean \quad MAD \quad {MSE}}$ &  $\dfrac{N=1000}{Mean \quad MAD \quad {MSE}}$
   \vspace{12pt}\\
  
  \hspace{5pt}$\hat{\alpha}$ & 0.7030\quad0.0327\quad0.0017&
  0.7022\quad0.0315\quad0.0015 &
  0.7008\quad0.0307\quad0.0015
 \vspace{12pt}\\
  
 \vspace{5pt} $\hat{\theta}$ & 0.4033\quad0.0139\quad0.000295 &  
 0.4009\quad0.0136\quad0.000249& 0.4007\quad0.0130\quad0.000266\\
\end{tabular}
\\\\
Table $4$
Mean estimates of and dispersions from the true parameter for a simulated InG subordinator data with ($\alpha$, $p$)=($0.9$, $0.7$).\\\\

\begin{tabular}{c c c c}
  &$\dfrac{N=100}{Mean \quad MAD \quad {MSE}}$ &  $\dfrac{N=500}{Mean \quad MAD \quad {MSE}}$ &  $\dfrac{N=1000}{Mean \quad MAD \quad {MSE}}$\vspace{12pt}\\
  $\hat{\alpha}$ & 0.9319\quad0.0394\quad0.0017&
  0.9216\quad0.0225\quad0.0006
  & 0.9071\quad0.0214\quad0.0005 
  \vspace{12pt}\\
  $\hat{p}$ & 0.6644\quad0.0446\quad0.0021 &  0.6753\quad0.0257\quad0.0008 & 0.6922\quad 0.0250\quad0.0007
  \vspace{12pt}\\
\end{tabular}

Next, we use the asymptotic behavior of fractional moment of InG and InG-$\epsilon$ subordinators  to estimate $\alpha$ 
and  compare them using MAD and MSE. 
 Recall we have asymptotic behaviour of $p^{th}$ fractional moment of InG subordinator (see \cite{ricciuti}) as
\begin{equation}\label{asym beh of InG}\mathbb{E}S_{\alpha}^{p}(t)\simeq\dfrac{\Gamma(1-\frac{p}{\alpha})}{\Gamma(1-p)}t^{\frac{p}{\alpha}},\qquad t\rightarrow\infty,\end{equation}
where $p\in (0,1]$ and $p\leq \alpha$.

  Here we generate $n=10$ samples of InG  and InG-$\epsilon$ subordinator using their algorithm  with sample size $N=50$, $100$ and $250$  and estimate $\alpha$ using \eqref{asym beh of InG}. These estimated values are shown along with their MAD and MSE. It can be observed from Table $5$ and $6$ that estimated values of estimator are very close to true value with hardly one percent of variation. \\\\
Table $5$
Mean estimates of and dispersions from the true parameter for a simulated InG subordinator data for parameter $\alpha$.\\

\begin{tabular}{c c c c c}
  &$\dfrac{N=50}{Mean \quad MAD  \quad {MSE}}$ &  $\dfrac{N=100}{Mean \quad MAD \quad {MSE}}$ &  $\dfrac{N=250}{Mean \quad MAD \quad {MSE}}$
   \vspace{12pt}\\
  
 $\alpha=0.3$ \hspace{5pt}& 
  0.3086\quad0.0183\quad0.000059 &
  0.3060\quad0.0176\quad0.00047
  & 0.3022\quad0.0207\quad0.000711 
 \vspace{12pt}\\

 $\alpha=0.5$ \hspace{5pt}&  0.5038\quad0.0061\quad0.000047 & 0.5019 \quad 0.0063 \quad 0.000066 &

  0.5000\quad0.0072\quad0.000087
 \vspace{12pt}\\
 
  $\alpha=0.7$ \hspace{5pt}& 
  0.7059\quad0.0126\quad0.00033 &
  0.7037\quad0.0120\quad0.0002392
  &0.7019\quad0.0127\quad0.000386
 \vspace{12pt}\\
  $\alpha=0.9$ \hspace{5pt}
& 
  0.9086\quad0.0071\quad0.000107
  & 0.9064\quad0.0114\quad0.00024 &0.9054\quad0.0094\quad0.00016
  \vspace{12pt}\\\\

\end{tabular}

Table $6$
Mean estimates of and dispersions from the true parameter for a simulated InG-$\epsilon$ subordinator data for parameter $\alpha$.\\

\begin{tabular}{c c c c c}
  &$\dfrac{N=50}{Mean \quad MAD  \quad {MSE}}$ &  $\dfrac{N=100}{Mean \quad MAD \quad {MSE}}$ &  $\dfrac{N=250}{Mean \quad MAD \quad {MSE}}$
   \vspace{12pt}\\
  
 $\alpha=0.3$ \hspace{5pt}& 
  0.3072\quad0.0163\quad0.00007 &
  0.3053\quad0.0158\quad0.00053
  & 0.3019\quad0.0127\quad0.00041 
 \vspace{12pt}\\

 $\alpha=0.5$ \hspace{5pt}&  0.5065\quad0.0071\quad0.00008 & 0.5027 \quad 0.0053 \quad 0.00006 &

  0.5007\quad0.0042\quad0.000062
 \vspace{12pt}\\
 
  $\alpha=0.7$ \hspace{5pt}& 
  0.7067\quad0.0152\quad0.00072 &
  0.7053\quad0.0140\quad0.00053
  &0.7015\quad0.0134\quad0.00029
 \vspace{12pt}\\
  $\alpha=0.9$ \hspace{5pt}
& 
  0.9096\quad0.0086\quad0.00026
  & 0.9056\quad0.0078\quad0.00021 &0.9042\quad0.0080\quad0.00012
  \vspace{12pt}\\\\

\end{tabular}

\paragraph{$\mathbf{Acknowledgment}$}  First author would like to acknowledge the Centre for Mathematical \& Financial Computing 
and the DST-FIST grant for the infrastructure support for the computing 
lab facility under the scheme FIST (File No: SR/FST/MS-I/2018/24) 
at the LNMIIT, Jaipur.
\bibliographystyle{abbrv}
\bibliography{refs}

@article{10.1111/j.2517-6161.1971.tb00884.x,
    author = {Macdonald, P. D. M.},
    title = {Comments and Queries Comment on “An Estimation Procedure for Mixtures of Distributions” by Choi and Bulgren},
    journal = {Journal of the Royal Statistical Society: Series B (Methodological)},
    volume = {33},
    number = {2},
    pages = {326-329},
    year = {2018},
    month = {12},
    issn = {0035-9246},
    doi = {10.1111/j.2517-6161.1971.tb00884.x},
    url = {https://doi.org/10.1111/j.2517-6161.1971.tb00884.x},
    eprint = {https://academic.oup.com/jrsssb/article-pdf/33/2/326/49098858/jrsssb\_33\_2\_326.pdf},
}

@Article{ricciuti,
AUTHOR = {Beghin, Luisa and Ricciuti, Costantino},
TITLE = {Lévy Processes Linked to the Lower-Incomplete Gamma Function},
JOURNAL = {Fractal and Fractional},
VOLUME = {5},
YEAR = {2021},
NUMBER = {3},
ARTICLE-NUMBER = {72},
URL = {https://www.mdpi.com/2504-3110/5/3/72},
ISSN = {2504-3110},
DOI = {10.3390/fractalfract5030072}
}

@Preamble{
"\def\cprime{$'$} "
}

@article{Murphy2015BayesianPE,
  title={Bayesian parameter estimation of Jump-Langevin systems for trend following in finance},
  author={James K. Murphy and Simon J. Godsill},
  journal={2015 IEEE International Conference on Acoustics, Speech and Signal Processing (ICASSP)},
  year={2015},
  pages={4125-4129},
  url={https://api.semanticscholar.org/CorpusID:14923625}
}

@inproceedings{Chalana1997ParameterEI,
  title={Parameter estimation in deformable models using Markov chain Monte Carlo},
  author={Vikram Chalana and David R. Haynor and Paul D. Sampson and Yongmin Kim},
  booktitle={Medical Imaging},
  year={1997},
  url={https://api.semanticscholar.org/CorpusID:62614560}
}

@article{Bhuvaneswari2019MachineLP,
  title={Machine Learning Parameter Estimation in a Smart-City Paradigm for the Medical Field},
  author={M. Bhuvaneswari and G. Naveen Balaji and Fadi M. Al-turjman},
  journal={Smart Cities Performability, Cognition, \& Security},
  year={2019},
  url={https://api.semanticscholar.org/CorpusID:198342701}
}

@article{babulal,
author = {Babulal, M. S. and Gauttam, S. K. and Maheshwari, A.},
title = {Lévy Processes with Jumps Governed by Lower Incomplete Gamma Subordinator and Its Variations},
journal = {Theory of Probability \& Its Applications},
volume = {70},
number = {1},
pages = {73-91},
year = {2025},
doi = {10.1137/S0040585X97T992240},

URL = { 
    
        https://doi.org/10.1137/S0040585X97T992240
    
    

},
eprint = { 
    
        https://doi.org/10.1137/S0040585X97T992240
    
    

}
}

@article{ALJEDDANI2023221,
title = {Parameter estimation of a model using maximum likelihood function and Bayesian analysis through moment of order statistics},
journal = {Alexandria Engineering Journal},
volume = {75},
pages = {221-232},
year = {2023},
issn = {1110-0168},
doi = {https://doi.org/10.1016/j.aej.2023.05.079},
url = {https://www.sciencedirect.com/science/article/pii/S1110016823004465},
author = {Sadiah M.A. Aljeddani and M.A. Mohammed}
}

@article{Cahoy2012SimulationAE,
  title={Simulation and Estimation for the Fractional Yule Process},
  author={Dexter O. Cahoy and Federico Polito},
  journal={Methodology and Computing in Applied Probability},
  year={2012},
  volume={14},
  pages={383-403},
  url={https://api.semanticscholar.org/CorpusID:59432161}
}

@article{garch,
    author = {Haug, S. and Klüppelberg, C. and Lindner, A. and Zapp, M.},
    title = {Method of moment estimation in the COGARCH(1,1) model},
    journal = {The Econometrics Journal},
    volume = {10},
    number = {2},
    pages = {320-341},
    year = {2007},
    month = {06},
    issn = {1368-4221},
    doi = {10.1111/j.1368-423X.2007.00210.x},
    url = {https://doi.org/10.1111/j.1368-423X.2007.00210.x},
    eprint = {https://academic.oup.com/ectj/article-pdf/10/2/320/25780888/ectj0320.pdf},
}

@article{articlefbpdpPE,
author = {Cahoy, Dexter and Polito, Federico},
year = {2012},
month = {03},
pages = {1-12},
title = {Parameter estimation for fractional birth and fractional death processes},
volume = {24},
journal = {Statistics and Computing},
doi = {10.1007/s11222-012-9365-1}
}

@article{MLEitoP,
 ISSN = {02664666, 14694360},
 URL = {http://www.jstor.org/stable/3532294},
 author = {Andrew W. Lo},
 journal = {Econometric Theory},
 number = {2},
 pages = {231--247},
 publisher = {Cambridge University Press},
 title = {Maximum Likelihood Estimation of Generalized Itô Processes with Discretely Sampled Data},
 urldate = {2025-06-24},
 volume = {4},
 year = {1988}
}

@article {Cahoy2010,
	AUTHOR = {Cahoy, Dexter O. and Uchaikin, Vladimir V. and Woyczynski,
	Wojbor A.},
	TITLE = {Parameter estimation for fractional {P}oisson processes},
	JOURNAL = {J. Statist. Plann. Inference},
	FJOURNAL = {Journal of Statistical Planning and Inference},
	VOLUME = {140},
	YEAR = {2010},
	NUMBER = {11},
	PAGES = {3106--3120},
	ISSN = {0378-3758},
	CODEN = {JSPIDN},
	MRCLASS = {62F10 (60G18 60J27 62E20)},
	MRNUMBER = {2659841 (2011h:62059)},
	DOI = {10.1016/j.jspi.2010.04.016},
	URL = {http://dx.doi.org/10.1016/j.jspi.2010.04.016},
}

@inbook{inbook,
author = {Saleh, A.K.Md.Ehsanes and Rohatgi, Vijay},
year = {2015},
month = {09},
pages = {},
title = {An Introduction to Probability and Statistics}
}
\end{document}